\documentclass[12pt, twoside, leqno]{article}

\usepackage{amsmath,amsthm}
\usepackage{amssymb}

\usepackage{enumerate}

\usepackage{graphicx}

\newtheorem{thm}{Theorem}[section]
\newtheorem{cor}[thm]{Corollary}
\newtheorem{lem}[thm]{Lemma}
\newtheorem{prop}[thm]{Proposition}

\theoremstyle{definition}

\newtheorem{rem}[thm]{Remark}

\numberwithin{equation}{section}

\usepackage{pgf,tikz}
\usetikzlibrary{arrows}

\begin{document}

\baselineskip=12.1pt

%%%%%%%%%%%%%%%%

\title{ The distance spectrum of the line graph of the crown graph}

\author{S.Morteza Mirafzal\\
Department of Mathematics \\
  Lorestan University, Khorramabad, Iran\\
E-mail: mirafzal.m@lu.ac.ir\\
E-mail: smortezamirafzal@yahoo.com
}

\date{}

\maketitle

\renewcommand{\thefootnote}{}


 % \footnote{\emph{*Corresponding author.}}
\footnote{\emph{Date}: }
\renewcommand{\thefootnote}{\arabic{footnote}}

\setcounter{footnote}{0}

%----------additions

\begin{abstract} 
The distance eigenvalues of a connected graph $G$ are the eigenvalues of its distance matrix
$D(G)$.  A graph is called distance integral if all of its
distance eigenvalues are integers. Let $n \geq 3$ be an integer.  The crown graph $Cr(n)$ is a graph obtained from the complete bipartite graph $K_{n,n}$ by removing a perfect matching. Let $L(Cr(n))$ denote the line graph of the crown graph $Cr(n)$. Using the equitable partition method,  the set of  distinct  distance eigenvalues of the graph $L(Cr(n))$ has been determined which shows  that this graph is distance integral [S.Morteza   Mirafzal, The line graph of the crown graph is distance integral,  Linear and Multilinear Algebra 71, no. 4 (2023): 662-672]. The     distance spectrum of the graph $L(Cr(n))$  has not  been found yet. In this paper, having the set of distance eigenvalues of $L(Cr(n))$ in the hand,  we determine the distance spectrum of this graph. 

\end{abstract}

\maketitle
\section{ Introduction}
 In this paper, a graph $G=(V,E)$ is
considered as an undirected simple graph where $V=V(G)=\{v_1,\dots,v_n\}$ is the vertex-set
and $E=E(G)=\{e_1,\dots,e_m\}$ is the edge-set. For all the terminology and notation
not defined here, we follow [11,12,16].\\
Let $G=(V,E)$ be a graph and  $A=A(G)$ be an adjacency matrix of $G$. The $characteristic$  $polynomial$  of $G$ is defined as $P(G; x)=P(x) = |xI-A|=det(xI-A)$, where $I$ is the identity matrix. 
Since $A$ is a real
symmetric matrix, the  characteristic polynomial $P(x)$ has real zeros. Every zero of the polynomial $P(x)$ is called an  $eigenvalue$  of the graph $G$. If $\lambda$ is a zero of $P(x)$, then the algebraic multiplicity of $\lambda$ is the multiplicity of it as a root of $P(x)$. 
The geometric   multiplicity of $\lambda$ is the dimension of its eigenspace.  Since $A$ is a real symmetric matrix, then these multiplicities are the same and this common number is called the $multiplicity$ of $\lambda$.
  A graph
is called $integral$   if all  of its eigenvalues are integers. The study of integral graphs  was initiated by Harary
and Schwenk in 1974 (see [17]). A survey of papers up to 2002 has been  appeared in [8], but
more than a hundred new studies on integral graphs have been published in the last
twenty three years. \newline 
The $distance$ between the vertices  $v_i$ and $v_j$, denoted by $d(v_i,v_j)$  is defined as the length of a shortest path between $v_i$ and $v_j$. The $ distance$ $matrix$  of $G$ denoted by $D(G)=D$ is the $n\times n$ matrix in which the rows and columns are indexed by the vertex-set whose $(i,j)$-entry is equal to $d(v_i,v_j)$  for $1\leq i,j\leq n$. The characteristic polynomial of $D(G)$ is defined $P_D(x)=P_{D(G)}(x)=Det(xI-D(G))$,  where $I$ is the  $n\times n$ identity matrix. It is called the $distance$  $characteristic$  $polynomial$  of $G$. Since $D(G)$ is a real symmetric matrix all of its eigenvalues  called $distance$   $eigenvalues$ of $G$, are real. The spectrum of $D(G)=D$  is denoted by
 $ Spec(D)=\{\lambda_1,\lambda_2,\dots,\lambda_n\}$ and indexed such that  $\lambda_1\geq\lambda_2\dots\geq \lambda_n$,  is called the $distance$  $spectrum$  of   $G$.  If the eigenvalues of $D$ are ordered by
$ \lambda_1 > \lambda_2 > \dots> \lambda_r  $, and their multiplicities are $ m_1, m_2, \dots,m_r $,  respectively,
 then we write,  \newline \newline
 \centerline{$ Spec( D) $   =  ${ \lambda_1,\lambda_2,\dots, \lambda_r } \choose{ m_1, m_2, \dots,m_r } $ 
 or  $Spec( D)  $ = $ \{  \lambda_1^{m_1},  \lambda_2^{m_2}, \dots, \lambda_r^{m_r}   \}.$ }  \newline \newline
  Since $D$ is an irreducible, non-negative, real and symmetric matrix,  from matrix theory it follows that  $\lambda_1$ is a simple eigenvalue and satisfies $\lambda_1\geq|\lambda_i|$, for $i=2,3,\dots,n$, and there exists a positive eigenvector corresponding to $\lambda_1$ [11,12,16]. The largest eigenvalue $\lambda_1$ is called the $distance$ $spectral$ $radius$ or $distance$  $index$ of the graph $G$.\\
The distance matrix and distance eigenvalues of graphs have been studied by researchers for many years. Some of the recent results and surveys concerning the subject include [1,2,3,4,5,6,7,9,10,13,14,15,18,19,31,33,34].
A graph $G$ is  $distance$  $integral$  (briefly, $D$-$integral$) if all the distance eigenvalues  of $G$ are integers. Although there are many
 papers that study distance spectrum of graphs and their applications, the $D$-integral graphs are studied only in a few number of papers (see [1,13,14,20,26,28,29,30,32,35]). \\
Let $n \geq 3$ be an integer.  A $crown$  $graph$  $Cr(n)$ is a graph obtained from the complete bipartite graph $K_{n,n}$ by removing a perfect matching.  The $bipartite$  $Kneser$  graph $H(n,k)$, $1 \leq k \leq n-1$, is a bipartite graph with the vertex-set consisting of all $k$-subsets and $(n-k)$-subsets of the set $[n]=\{1,2,3,\dots,n\}$, in which two vertices $v$ and $w$ are adjacent if and only if $v \subset w$ or $w \subset v$. Recently, this class of graphs has been studied from several  aspects [21,22,23,24,25]. 
It is easy to see that the crown graph $Cr(n)$ is isomorphic with the bipartite graph 
$H(n,1)$.   It is easy to check that the crown graph $Cr(n)$ is a vertex and edge-transitive graph of order $2n$ and regularity $n-1$ with diameter 3. In fact, $Cr(n)$ is a distance transitive graph [21,27]\
 It has been shown that the graph $Cr(n)$ is a distance integral graph [20,28,30]. 
 Let $L(Cr(n))$ denote the line graph of the crown graph $Cr(n)$.  It is not hard to see that the graph $L(Cr(n))$ is a vertex-transitive graph of order $n(n-1)$ and regularity $2(n-1)-2=2n-4$ with diameter 3.\
From the various interesting properties  of the graph $L(Cr(n))$, we are  interested in its distance spectrum.  Determining the set of distinct distance eigenvalues of the graph $L(Cr(n))$,   it has been proved that this graph is distance integral [26]. Up to our knowledge, the     distance spectrum of the graph $L(Cr(n))$  has not  been found yet.  In this paper, having the set of distinct distance eigenvalues of $L(Cr(n))$ in the hand,  we determine the distance spectrum of this graph. 
\section{Preliminaries }
The set of all permutations of a set $V$ is denoted by $Sym(V)$. A permutation group on $V$ is a subgroup of $Sym(V)$. If $G=(V,E) $ is a graph, 
 then we can view each automorphism as a
permutation of
$ V,$  and so $Aut (G)$ is a permutation group.
A permutation representation of a group $\Gamma$ is a homomorphism from $\Gamma$
into $Sym(V)$ for some set  $V.$ A permutation representation is also referred
to as an $action$  of $\Gamma $ on the set $V, $ in which case we say that $\Gamma$ acts on $V.$ A permutation group $\Gamma$ on  $V$  is $transitive$ if given any two elements $x$ and
$y$ from $V$ there is an element $g \in \Gamma$ such that $x^g=y.$  For each $v\in V,$  the set
$v^{\Gamma}=\{ v^g \ | \ g \in \Gamma  \}$ is called an $orbit$  of  $\Gamma$.  
 It is easy to see that if $\Gamma$ acts on $V$, then $\Gamma$ is
transitive  on $V$ (or $\Gamma$ acts $transitively$ on $V$), when there is just
one orbit. It is easy to see that the set of orbits of $\Gamma$ on $V$ is a partition of the set $V$.\\
 A graph $G=(V,E)$ is called $vertex$-$transitive$ if  $Aut(G)$
acts transitively on $V.$  We say that $G$ is $edge$-$transitive$ if the group $Aut(G)$ acts transitively  on the edge set $E$, namely, for any $\{x, y\} ,   \{v, w\} \in E(G)$, there is some $a$ in $Aut(G)$,  such that $a(\{x, y\}) = \{v, w\}$.   We say that $G$ is $distance$-$transitive$ if  for all vertices $u, v, x, y$ of $G$ such that $d(u, v)=d(x, y)$, where $d(u, v)$ denotes the distance between the vertices $u$ and $v$  in $G$,  there is an automorphism $\alpha$ in $Aut(G)$ such that  $\alpha(u)=x$ and $\alpha(v)=y.$ \\
 Let  $\alpha$ be a permutation of a set $X=\{ x_1, \dots,x_n \}$. This permutation can be represented by a permutation matrix $P=(p_{ij})$, where $p_{ij}=1$ if $\alpha(v_j)=v_i$, and $p_{ij}=0 $ otherwise. In the sequel, we need the following fact.
\begin{prop}  $[11]$
Let $A$ be the adjacency matrix of a graph $G=(V,E)$, and $\alpha$ a permutation of $V$. Then $\alpha$ is an automorphism of the graph $G$ if and only if $PA = AP$, where $P$ is the permutation matrix representing $\alpha$.

\end{prop}
If $G_1=(V_1,E_1),G_2=(V_2,E_2)  $ are  graphs, then their direct product is the graph
$  G_1 \times G_2 $
 with the  vertex-set $ \{( v_1,v_2) \  |  \  v_1 \in G_1,   v_2 \in G_2 \} $, and for
  which vertices $( v_1,v_2)$ and $ ( w_1,w_2)  $ are adjacent precisely if $ v_1$ is adjacent to $w_1$ in $G_1$ and $ v_2$ is adjacent to $w_2$ in $ G_2$. When $G_2=K_2$, the complete graph on two vertices,  then $G \times K_2$ is known as the $double$ $cover$ of the graph $G$. It is a fact that if the spectrum of $G$
 is $  \lambda_1,\lambda_2, \dots, \lambda_n $, then the spectrum of the double cover of it, that is, $G \times K_2$ is $  \lambda_1,\lambda_2, \dots, \lambda_n,  -\lambda_1,-\lambda_2, \dots, -\lambda_n$ [12]. In other words, when $\lambda$ is an eigenvalue of the graph $G$, then $\lambda$ and $-\lambda$ are eigenvalues of $G \times K_2. $\\
 It is not difficult to show that the crown graph  $Cr(n)$ is isomorphic with the double  cover of the complete graph $K_n$ [22,25].  We know that the spectrum of $K_n$ is\\ $ \{ {(n-1)}^1, {(-1)}^{n-1} \}$. Thus  the spectrum of the
  crown graph  is
  $$\{ {(n-1)}^1, {1}^{n-1},{(-1)}^{n-1}, {(-n+1)}^1 \}.$$ 
\section{Main results} 
Let $n \geq 3$ be an integer and $[n]= \{1,2,\dots,n \}$. Let $X= \{x_1,x_2,\dots,x_n  \}$ be an $n$-set disjoint from $[n]$.
 It is easy to see that the crown graph $Cr(n)$ is a graph with the vertex-set $[n] \cup X$ and the edge-set $E_0=\{e_{ij}=\{i,x_j \} \ | \ i,j \in [n], i \neq j  \}.$   Thus, $L(Cr(n))$, the line graph of $Cr(n)$,  is a graph with the vertex-set $E_0$ in which two vertices $e_{ij}$ and $e_{rs} $ are adjacent if and only if $i=r$ or $j=s$.\ %
  Let $V=\{(i,j) \ | \ i,j \in [n], i \neq j  \}$. Let $G$ be a graph with the vertex-set $V$ in which two vertices $(i,j)$ and $(r,s)$ are adjacent if and only if $i=r$ or $j=s$.
  It is easy to check that the graph $G$ is isomorphic with the graph $L(Cr(n))$.
  Hence,  in the sequel we work on the graph $G$ and call it the line graph of  the crown graph $Cr(n)$ and denote it by $L(Cr(n))$.\
It is easy to see that $L(Cr(3))$ is the cycle  graph $C_6$, which its structure is known and its line graph is again $C_6$. Hence in the rest of the paper we assume that $n\geq 4.$  It is easy to check that two non adjacent  vertices $(i,j)$ and $(r,s)$  are at distance 2 from each other whenever $i=s$ or $j=r$ or 
  $\{i,j\} \cap \{r,s\}=\varnothing$. Moreover, vertices  $(i,j)$ and $(j,i)$ are at distance 3 from each other ($P: (i,j),(x,j),(x,i),(j,i)$ is a shortest path between the vertices $(i,j)$ and $(j,i))$.   Thus the diameter of the graph $L(Cr(n))$ is 3. \
    Note that   the crown graph $Cr(n)$ is a regular graph and its adjacency spectrum is known,  
   hence the adjacency spectrum of its line graph, that is, $L(Cr(n))$ is known [11,12,16].
  \begin{rem} Although the crown graph $Cr(n)$ is a distance-transitive graph (and consequently it is distance-regular [11,16,21]), it is easy to check that the graph $L(Cr(n))$ is not distance-regular. Hence we can not use the theory of distance-regular graphs for determining the set of distance eigenvalues or spectrum.
  \end{rem} 
   Figure 1  shows the graph $L(Cr(4))$. Note that  the vertex $(i,j)$
 is denoted by $ij$ in this figure. \\
\definecolor{qqqqff}{rgb}{0.,0.,1.}
\begin{tikzpicture}[line cap=round,line join=round,>=triangle 45,x=0.6384255715985469cm,y=.63cm]
\clip(3.371,-7.1306) rectangle (22.0292,3.0576);
\draw (9.8566,2.6462) node[anchor=north west] {41};
\draw (14.7692,2.5736) node[anchor=north west] {42};
\draw (16.8746,-0.0642) node[anchor=north west] {12};
\draw (16.8746,-2.5568) node[anchor=north west] {13};
\draw (15.8098,-5.2188) node[anchor=north west] {14};
\draw (9.9292,-5.4366) node[anchor=north west] {34};
\draw (7.5092,-2.7746) node[anchor=north west] {32};
\draw (6.9914,-0.3788) node[anchor=north west] {31};
\draw (9.1114,0.3714) node[anchor=north west] {21};
\draw (12.8332,-4.5896) node[anchor=north west] {24};
\draw (14.8902,-2.0002) node[anchor=north west] {23};
\draw (10.42,1.8)-- (14.7692,1.8476);
\draw (14.7692,1.8476)-- (16.52,-0.7);
\draw (16.52,-0.7)-- (16.5,-2.94);
\draw (16.5,-2.94)-- (15.08,-5.24);
\draw (15.08,-5.24)-- (10.54,-5.22);
\draw (10.54,-5.22)-- (8.22,-2.8);
\draw (8.22,-2.8)-- (8.16,-0.7);
\draw (8.16,-0.7)-- (10.42,1.8);
\draw (10.42,1.8)-- (10.32,-0.12);
\draw (10.32,-0.12)-- (8.16,-0.7);
\draw (12.86,-4.34)-- (10.54,-5.22);
\draw (12.86,-4.34)-- (15.08,-5.24);
\draw (14.04,-0.66)-- (14.48,-2.46);
\draw (14.7692,1.8476)-- (14.04,-0.66);
\draw (14.2318,0.1326) node[anchor=north west] {43};
\draw (16.5,-2.94)-- (14.04,-0.66);
\draw (14.48,-2.46)-- (16.5,-2.94);
\draw (10.32,-0.12)-- (14.48,-2.46);
\draw (16.52,-0.7)-- (15.08,-5.24);
\draw (10.32,-0.12)-- (12.86,-4.34);
\draw (16.52,-0.7)-- (8.22,-2.8);
\draw (14.48,-2.46)-- (12.86,-4.34);
\draw (10.42,1.8)-- (14.04,-0.66);
\draw (14.7692,1.8476)-- (8.22,-2.8);
\draw (8.16,-0.7)-- (10.54,-5.22);
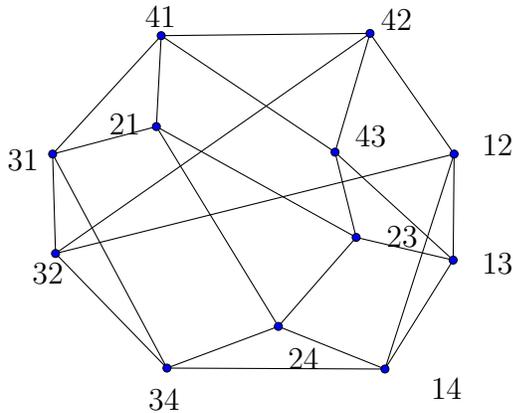
\draw (8.4328,-6.078) node[anchor=north west] {Figure 1. The graph $L(Cr(4))$};
\begin{scriptsize}
\draw [fill=qqqqff] (10.42,1.8) circle (1.5pt);
\draw [fill=qqqqff] (14.7692,1.8476) circle (1.5pt);
\draw [fill=qqqqff] (16.52,-0.7) circle (1.5pt);
\draw [fill=qqqqff] (16.5,-2.94) circle (1.5pt);
\draw [fill=qqqqff] (15.08,-5.24) circle (1.5pt);
\draw [fill=qqqqff] (10.54,-5.22) circle (1.5pt);
\draw [fill=qqqqff] (8.22,-2.8) circle (1.5pt);
\draw [fill=qqqqff] (8.16,-0.7) circle (1.5pt);
\draw [fill=qqqqff] (10.32,-0.12) circle (1.5pt);
\draw [fill=qqqqff] (12.86,-4.34) circle (1.5pt);
\draw [fill=qqqqff] (14.04,-0.66) circle (1.5pt);
\draw [fill=qqqqff] (14.48,-2.46) circle (1.5pt);
\end{scriptsize}
\end{tikzpicture}\\
Since the graph $Cr(n)$ is distance-transitive, hence it is edge-transitive. Thus the graph $L(Cr(n))$ is a vertex-transitive graph. \\
Let $G=(V,E)$ be $k$-regular graph. It is easy to see that  the graph $L(G)$, the line graph of $G$, is a $(2k-2)$-regular graph. There is a close relationship between the spectrum of $G$ ant its line graph $L(G)$. 
\begin{prop} $[11]$
Let $G=(V,E)$ be a $k$-regular graph with $n$ vertices and $ m =\frac{1}{2}nk $ edges. If 
$$ Spec(G)= \{ {k}^{1},{\lambda_2}^{m_2}, \dots, {\lambda_t}^{m_t}  \}, $$
then for the spectrum of $L(G)$ we have,
$$ Spec(L(G))= \{ {(2k-2)}^1, {(k-2+\lambda_2)}^{m_2}, \dots, {(k-2+\lambda_t)}^{m_t}, {(-2)}^{m-n} \}.$$

\end{prop}
  
  In proving our main result, we need the following important result.      
 \begin{thm} $[26]$
Let $n>3$ be an integer. Then the line graph of the crown graph, that is, the graph $L(Cr(n))$,  is a distance  integral graph with the set of  distance eigenvalues,  $$S=\{-n-1, -n+3, - 1, 1,  2n^2-4n+3\}. $$
 
 \end{thm}                      
 In the sequel, $J$ is the all 1 matrix of appropriate size,  and $I $ is the identity matrix.\\ 
Let $G=(V,E)$, $V=\{v_1, \dots v_n  \}$,  be a connected graph  with diameter $d$.
For every integer $i$,  $0\leq i \leq d$, the distance-$i$ matrix $A_i$ of $G$ is defined as a matrix whose
rows and columns are indexed by the vertex-set of $G$ and the entries are defined as, 
\begin{center}
$A_{i}(v_r,v_s)=\begin{cases}1 & \ if \ d(v_r,v_s)=i \\ 0 & \ otherwise.  \end{cases} $
\end{center}
Then $A_0 = I$, $A_1$ is the usual adjacency matrix $A$ of $G$. Note   that $$A_0 + A_1 + \dots + A_d = J,$$
  Now it is clear that if $D(G)$ is the distance matrix of $G$, then 
   $$D(G)= A_1 + 2A_2  + 3A_3 + \dots + dA_d. \ \ \ \ \  (1)$$ 
\begin{lem}
Let $G=(V,E)$ be a graph with the adjacency matrix $A$ and distance matrix $D$.  
Let  $A_i$  be the  distance-$i$ matrix  of $G$. If the diameter of $G$ is $3$, then we have 
$$ D=-A+2J-2I+A_3  $$
\end{lem}

\begin{proof}

Since  the diameter of the graph $G$ is 3, hence by (1),  we have $D=A+2A_2+3A_3$. Let $\overline{G}$ be the complement of the graph $G$ with the adjacency matrix $\overline{A}$. It is clear that $\overline{A}+A=J-I$, thus $\overline{A}=J-I-A$. On the hand, it is easy to see that $\overline{A}=A_2+A_3$. Hence we have $ A_2 +A_3=J-I-A, $ and thus $A_2=J-I-A-A_3.$ We now deduce  that
 $$D=A+2(J-I-A-A_3)+3A_3= -A+2J-2I+A_3.  \ \ \ (2) $$

\end{proof}
We know that the diameter of the graph $G=L(Cr(n))$ is 3, hence by Lemma 3.4, we have the following corollary.
\begin{cor}
Let $n \geq 3$ be an integer and $G=L(Cr(n))$ be the line graph of the crown graph $Cr(n)$ with the adjacency matrix $A=A(G)$ and distance matrix $D=D(G)$. Then we have $$D=-A+2J-2I+A_3,$$  where $  A_3 $  is the distance-$3$ matrix of $G$, $J$ is the matrix with all entry 1 and $I$ is the identity matrix, 
\end{cor}

In the sequel, we need the spectrum of $A_3$ for the graph $G=L(Cr(n))$. By the following result we determine  the spectrum of $A_3.$

\begin{lem}
Let $A$ be the adjacency matrix and  $A_3$ be the distance-$3$ matrix of the graph $G=L(Cr(n))$. Let $m=n(n-1)$ and $k=\frac{m}{2}$. Then \\
  $$AA_3=A_3A,$$
   and for the spectrum of $A_3$ we have, 
 $$ Spec(A_3)   =   \{  1^{k},  {-1}^{k} \}.$$ 

\end{lem}

\begin{proof}
Let $v=(i,j)$ be a vertex in the graph $G=L(Cr(n))$. It is easy to check that the vertex $v^c=(j,i)$ is the unique vertex in $G$ such that $d(v,v^c)=3$.
Let $V=V(G)$ be the vertex-set  of $G$. Consider the following mapping $\alpha, $ 
$$\alpha: V \rightarrow V, \ \ \alpha(v)=v^c, \ \  \forall v\in V. $$ 
It  is not difficult to check that $\alpha$  is an automorphism of the graph $G$. 
 We now can see since $A$ is the  adjacency matrix for $G$ and $A_3$ is  the  permutation matrix  of the automorphism  $\alpha$, then by  Proposition 2.1, we have $AA_3=A_3A.$ \\      Also, since ${\alpha}^2= \epsilon$, where $\epsilon$ is the identity mapping, hence ${A_3}^2=I$. Thus, if $\lambda$ is an eigenvalue of $A_3$,  
then ${\lambda}^2=1$. Hence $\lambda \in \{ 1,-1 \}$. Note that  each diagonal entry of $A_3$ is $0$. Thus, $trace(A_3)=0.$ It follows that if the multiplicity of the eigenvalue of 1 is $a$ and the eigenvalue of $-1$ is $b$, then $a-b=0$, hence $a=b$. Noting that $a+b=m$, we conclude that $a=b=\frac{1}{2}m$.

\end{proof}
We now proceed to prove and obtain the main result of our paper. 
 \begin{thm}
 Let $n \geq 3$ be an integer and $G=L(Cr(n))$ be the line graph of the crown graph $Cr(n)$ with the adjacency matrix $A=A(G)$ and distance matrix $D=D(G)$. Let $m=n(n-1) = |E(Cr(n))|$=$|V(L(Cr(n)))|$. Then for the spectrum of $D$ we have, 
 $$ Spec(D)=\{ {(2n^2-4n+3)}^1, 1^a, {(-1)}^{b},{(-n+3)}^{n-1},{(-n-1)}^{n-1} \},   $$
 where $a=\frac{1}{2}(n^2-3n+1)-\frac{1}{2}=\frac{1}{2}(n^2-3n)$ and\\ $b=\frac{1}{2}(n^2-3n+1)+\frac{1}{2}=\frac{1}{2}(n^2-3n+2).$
  \end{thm} 

\begin{proof}
Let $A=A_{m\times m}$ and $D=D_{m\times m}$ be the adjacency and distance matrix of the graph $G=L(Cr(n))$, respectively. By Lemma 3.4, we know that 
$$ D=-A+2J-2I+A_3,   $$ 
where $J=J_{m \times m} $, $ I= I_{m \times m}$ and $A_3$ are the matrices  with  all entry 1, the identity matrix and distance-$3$ matrix of $G$, respectively. Lemma 3.6, follows that 
  $${A_3}^2=I, \  A_3A=AA_3.$$
   Since the crown graph $Cr(n)$ is an $(n-1)$-regular graph, then $G=L(Cr(n))$ is a $(2n-4)$-regular graph. Thus we have $AJ=JA$ [11]. It is easy to see that $A_3J=JA_3 =J.$ We now deduce that the set,
     $$T= \{ A, J, I, A_3  \},   $$
  is a  commuting set  of matrices on $\mathbb{R}$,  the field  of real numbers. It is clear that each element of $T$ is a symmetric matrix. 
 Hence,  there is a  basis
  $$ B= \{e_1, \dots,e_m  \} $$
  of $ { \mathbb{R}}^m$ such that
 each $ e_i $ is an eigenvector for all of  the matrices in $T$ [12,16]. 
 Let  $j$ be a column of the matrix $J$,   Since the graph $G$ is a $(2n-4)$-regular graph, hence $Aj=(2n-4)j.$ We now have, \\
  $Dj=(-Aj)+2Jj-2Ij+(A_3)j=(-(2n-4)+2m-2+1)j$=$(-2n+4+2n^2-2n-1)j$=
$(2n^2-4n+3)j.$ \\
 Thus, $(2n^2-4n+3)$ is an eigenvalue of $D$.
We can assume that $e_1=j$  (or any element $e_r$ in $B$ such that $De_r=(2n^2-4n+3)e_r$. As we can check there is a unique element $e_r$ in $B$  such that $Je_r=me_r=n(n-1)e_r$. Also, for such an $e_r$,  we must have $Ae_r=(2n-4)e_r$.). In fact, since the rank of the  matrix $J$ is 1 and $Je_1=me_1$, we deduce that
  $$ Spec(J)= \{ m^1,0^{m-1} \}.  $$ 
The crown graph $Cr(n)$ is isomorphic with the graph $  K_n \times K_2$, where $K_k$ is the complete graph on $k$ vertices, hence
 $$ Spec(Cr(n))=\{ {(n-1)}^1, 1^{n-1}, {(-1)}^{n-1}, {(-n+1)}^1  \}.   $$ 
We  know that $G=Cr(n)$ is a $(n-1)$-regular graph. Now, noting that $n-3-n+1=-2$ and   $n(n-1)-(2n-1)=n^2-3n+1$,   Proposition 3.2, follows  that 
$$Spec(G)=Spec(L(Cr(n)))= \{ {(2n-4)}^1,  {(n-2)}^{n-1},  {(n-4)}^{n-1},  {(-2)}^{n^2-3n+1}   \}. \ \  (3)   $$ 
We recall that by Theorem 3.3, the set of distinct eigenvalues of the matrix $D$  is $$S=\{-n-1, -n+3, - 1, 1,  2n^2-4n+3\}. \ \ \ \ \ \ \ \ (4)$$ 
We now determine the multiplicity of each of the eigenvalue.\\
(i) Noting that $(A_3)j=j$  we have,  \\ $De_1=Dj=(-Aj)+2Jj-2Ij+(A_3)j=(-(2n-4)+2m-2+1)j$=$(-2n+4+2n^2-2n-1)j$=
$(2n^2-4n+3)j.$ \\ We know since $G$ is a regular graph, then the multiplicity of the corresponding eigenvalue of $j$ is 1. Hence the multiplicity of the distance eigenvalue 
$2n^2-4n+3$ is 1.\\
(ii) Let $e_i \in B$ is an eigenvector of $A$ corresponding to the eigenvalue $n-2.$ We now have,   
 $De_i=(-Ae_i)+2Je_i-2Ie_i+(A_3)e_i$=$(-n+2)e_i+0e_i-2e_i+A_3e_i.$  \\
We know by Lemma 3.6,  that  $$A_3e_i \in \{ -e_i,e_i \}.$$
If $A_3e_i=e_i$, then   $De_i=(-n+2)e_i+0e_i-2e_i+e_i=(-n+1)e_i$. \\
 Hence,  $(-n+1)$ is a distance eigenvalue of $G$, that is,  $-n+1 \in S$, which is a contradiction. Hence we must have $A_3e_i=-e_i$. We now have, \\
  $De_i=(-n+2)e_i+0e_i-2e_i-e_i=(-n-1)e_i$. \\
   We now deduce by (3) that  
   if the multiplicity of the distance eigenvalue  $(-n-1)$ is  $a_2$, then $a_2 \geq (n-1).$\\
(iii) Let $e_i \in B$ is an eigenvector of $A$ corresponding to the eigenvalue $n-4.$ We now have, \\
 $De_i=(-Ae_i)+2Je_i-2Ie_i+(A_3)e_i$=$(-n+4)e_i+0e_i-2e_i+A_3e_i.$\\
We know that,  $A_3e_i \in \{ -e_i,e_i \}.$
If $A_3e_i=-e_i$, then\\ $De_i=(-n+4)e_i+0e_i-2e_i-e_i=(-n+1)e_i$. \\ Thus,  $(-n+1)$ is a distance eigenvalue of $G$, that,  is $-n+1 \in S$, which is a contradiction. Hence we must have $A_3e_i=e_i$. We now have\\ $De_i=(-n+4)e_i+0e_i-2e_i+e_i=(-n+3)e_i$. \\
We now deduce by (3) that if the multiplicity of the distance eigenvalue  $(-n+3)$ is $a_3$, then $a_3 \geq (n-1).$\\
(iv) Let $e_i \in B$ is an eigenvector of $A$ corresponding to the
 eigenvalue $-2$, that is, $Ae_i=-2e.$ Hence we have \\
$De_i=(-Ae_i)+2Je_i-2Ie_i+(A_3)e_i$=$(2e_i)+0e_i-2e_i+A_3e_i = A_3e_i.$\\ If $A_3e_i=e_i$, then $De_i=e_i$, and if $A_3e_i=-e_i$, then $De_i=-e_i. $ Now it follows by (3) that if
the multiplicities of the distance eigenvalues  $1$ and $(-1)$ are $a$ and $b$, respectively, since  
$$ m-1=a_2+a_3+a+b \geq (n-1)+(n-1)+n^2-3n+1 =m-1,  $$
hence we must have $a_2=a_3=n-1.$\\
  Thus we have 
$$a+b = n(n-1)-(2n-2)-1 =n^2-3n+1. \ \ \ \ \ \ \\ \ (5)$$
 On the other hand, we know that  $trace(D)$   is the sum of all of its eigenvalues, hence from (4) we have\\ 
$0=trace(D)=(2n^2-4n+3)+a-b+(n-1)(-n+3)+(n-1)(-n-1)=$\\
$(2n^2-4n+3)+(-n^2+3n+n-3-n^2+1)+a-b=1+a-b.$
  Thus, we have
   $$a-b=-1 \ \ \ \ \ \ \ \ \ \ \ \ \ \ \ \ \ \ \ \ \ \ \ \ \  (6)$$ 
From (5) and (6) it follows that $2a=n^2-3n$, hence $a=\frac{1}{2}(n^2-3n)$ and $b=\frac{1}{2}(n^2-3n+2)$.  We now conclude the main result of the theorem, that is,
$$ Spec(D)=\{ {(2n^2-4n+3)}^1, 1^a, {(-1)}^{b},{(-n+3)}^{n-1},{(-n-1)}^{n-1} \},   $$
where $a=\frac{1}{2}(n^2-3n+1)-\frac{1}{2}=\frac{1}{2}(n^2-3n)$ and\\ $b=\frac{1}{2}(n^2-3n+1)+\frac{1}{2}=\frac{1}{2}(n^2-3n+2)$.
 
\end{proof}

\begin{rem}
Noting to the proof of Theorem 3.7,  we learn that the key point in the proof is having the set of distinct distance eigenvalues of the graph $L(Cr(n))$ in the hand which has been found in [26] by the equitable partition method.
\end{rem}

\section{ Acknowledgements}
The author is thankful to Oly Newman (University of Essex)  for a valuable comment on an earlier version of Theorem 3.7.

\end{document}